\documentclass[aps,preprint,amssymb,amsmath,eqsecnum,nofootinbib]{revtex4}
\usepackage{graphics}

\newtheorem{example}{Example}[section]
 \begin{document}
 \title{The non-adiabatic classical geometric phase and its bundle-theoretic interpretation}
 \author{Gavriel Segre}
 \homepage{http://www.gavrielsegre.com}
 \email{info@gavrielsegre.com}
 \date{6-4-2004}
 \maketitle
 \newpage
 \tableofcontents
 \newpage
 \section{Notation}
 \begin{tabular}{|c|c|}
   $ \sigma_{Borel} $ & Borel $ \sigma $-algebra of a topological
   space \\
   $ U_{T} $ & Koopman unitary of the dynamical system $ ( X  \, ,
   \, \sigma \, , \, \mu \, , \, T )$ \\
   P(M,G) & principal bundle with total-space P, base-space M and structure-group G \\
   ${\mathcal{S}} ({\mathcal{H}})$ & space of the normalized vectors in $ {\mathcal{H}} $  \\
   $P{\mathcal{H}}$ & space of the rays in $ {\mathcal{H}} $ \\
   $C_{p}(M)$ & space of the p-based loops in M  \\
   $ \omega_{Stiefel} $ & Stiefel connection on a universal bundle \\
   $ \tau_{\gamma}^{\omega} $ & holonomy of the loop $ \gamma $ w.r.t. the connection $ \omega $ \\  \hline
 \end{tabular}
 \newpage
 \section{Introduction}
 While the investigation of Berry phase has led to  \cite{Berry-90}, \cite{Bohm-93}, \cite{Mostafazadeh-02}  :
\begin{enumerate}
  \item a geometric interpretation of the phenomenon in terms of
  holonomies of a suitable connection over a suitable principal bundle over the parameters'-space
  (performed by Barry Simon)
  \item a nonadiabatic generalization of its (performed by Yakir
  Aharonov and Jeeva Anandan)
  \item a geometric interpretation of the nonadiabatic
  Aharonov-Anandan phase in terms of holonomies of the Stiefel
  connection on the U(1)-universal-bundle (performed by Yakir
  Aharonov and Jeeva Anandan themselves and refined and improved by Arno Bohm and
  Al\'{i} Mostafazadeh)
\end{enumerate}
only the first of these  conceptual step has been performed as to
classical Hannay phase \cite{Hannay-90} whose bundle-theoretic
interpretation has been shown by Ennio Gozzi and William D.
Thacker \cite{Gozzi-Thacker-87}.

Since Hannay's analysis applies to a very particular class of
dynamical systems (the integrable hamiltonian ones) I will follow
here a radically different approach consisting in:
\begin{enumerate}
  \item considering the holonomies of the Stiefel connection on
  the  U(1)-universal-bundle associated to the Koopman representation of
  an arbitrary classical dynamical system
  \item investigating their physical meanings as non adiabatic
  classical geometric phases
  \item obtaining the Hannay's result as a particular case
\end{enumerate}
\newpage
\section{Holonomies in the Koopman representation}

Given an arbitrary classical dynamical system $ ( X \, , \, \sigma
\, , \,  \mu \, ,  \, T ) $ (where hence $ ( X \, , \, \sigma \, ,
\, \mu ) $ is a classical probability space and  T is an
automorphism of its) \cite{Kornfeld-Sinai-00}, \cite{Benatti-93}
let us introduce the \textbf{Koopman unitary} $ U_{T} :
{\mathcal{H}} \, \mapsto \, {\mathcal{H}} $:
\begin{equation}
  ( U_{T} f ) (x) \; := \; ( f \circ T ) (x) \; \; f \in {\mathcal{H}}
\end{equation}
with:
\begin{equation}
  {\mathcal{H}} \; := \; L^{2} ( X , \mu)
\end{equation}
Let us now consider the \textbf{U(1)-universal-bundle} (cfr. e.g.
the mathematical appendix \cite{Mostafazadeh-93} of
\cite{Bohm-93} or \cite{Nakahara-95}) $ {\mathcal{S}}
({\mathcal{H}}  ) ( P({\mathcal{H}}) \, , \, U(1) ) $ whose
\textbf{base-space} (i.e. the \textbf{U(1)-classifying-space}) is
the space of all the rays in $ {\mathcal{H}} $, whose
\textbf{total space} $ {\mathcal{S}} ({\mathcal{H}} ) $ is the
unit sphere in $ {\mathcal{H}} $:
\begin{equation}
  {\mathcal{S}} ({\mathcal{H}}  ) \; := \; \{  | \psi >  \in
  {\mathcal{H}} \:  < \psi | \psi > \; = \; 1 \}
\end{equation}
and whose \textbf{structure group} is clearly U(1).

Endowed the U(1)-universal bundle $ {\mathcal{S}} ({\mathcal{H}}
) ( P {\mathcal{H}} \, , \, U(1) ) $ with the Stiefel connection
$ \omega_{Stiefel} $ let us denote by $ HOL_{Stiefel} ( | \psi > )
$ the \textbf{holonomy group in $ | \psi > $ } of $
\omega_{Stiefel} $:
\begin{equation}
  HOL_{Stiefel} ( | \psi > ) \; := \; \{ e^{i \theta} \in U(1) : \tau_{\gamma}^{\omega_{Stiefel}}
  ( | \psi > ) = e^{i \theta} | \psi > \, \}  \; \; \gamma \in C_{| \psi > < \psi|
  } (   P{\mathcal{H}}  )
\end{equation}
where $ C_{| \psi > < \psi|
  } (   P{\mathcal{H}}  ) $ is the space of the loops in $  P{\mathcal{H}} $  at $  | \psi > < \psi| $  while $ \tau_{\gamma}^{\omega_{Stiefel}} \, : \, \pi^{- 1} ( | \psi > < \psi| )
\mapsto \pi^{- 1} ( | \psi > < \psi| ) $ is the transformation of
the fibre in $ | \psi > < \psi| $ associated to a loop $ \gamma :
[ 0 , 1 ] \mapsto  P{\mathcal{H}} \: \gamma (0) = \gamma (1) = |
\psi > < \psi| $ in $ P{\mathcal{H}} $
\begin{example} \label{ex: uncoupled oscillators}
\end{example}
UNCOUPLED OSCILLATORS

Let us consider the hamiltonian dynamical system with
configuration space $ {\mathbb{R}}^{2} $  and hamiltonian  $ H :
T^{\star} {\mathbb{R}}^{2}  \, \mapsto \, {\mathbb{R}} $ given by:
\begin{equation}
  H( \vec{q} , \vec{p} ) \; := H_{1} + H_{2}
\end{equation}
\begin{equation}
  H_{i} ( \vec{q}_{i} , \vec{p}_{i} ) \; = \; \frac{1}{2} (
  \vec{p}_{i}^{2} + \omega_{i}^{2}  \vec{q}_{i}^{2} ) \; \; i =
  1,2
\end{equation}
The system is obviously integrable its \textbf{action variables}:
\begin{equation}
  \vec{I} \; := \; ( \frac{h_{1}}{\omega_{1}} \, , \,
  \frac{h_{2}}{\omega_{2}} )
\end{equation}
being constants of the motion:
\begin{equation}
 \frac{d \vec{I} }{ d t}  \; = \; 0
\end{equation}
while the \textbf{angle variables} $ \vec{\phi} $ evolve on the
2-torus $ T^{2} $ according to:
\begin{equation}
  T_{t}  ( \vec{\phi} ) \; = \;  \vec{\phi} + \vec{\omega}
  t
\end{equation}
where obviously $ \vec{\omega} := ( \omega_{1} \, , \, \omega_{2}
) $.

Since we have seen that our dynamical system may be seen as $ (
T^{2} , \sigma_{Borel} , d \mu ( \phi_{1} , \phi_{2} ) := \frac{d
\phi_{1} d \phi_{2}}{ (2 \pi)^{2} } \, , \, T_{t} ) $ the passage
to the Koopman representation involves the unitary operator $
U_{T_{t}} $ over $ {\mathcal{H}} := L^{2} ( T^{2} , d \mu (
\phi_{1} , \phi_{2} ) ) $ identified by its action on the basis:
\begin{equation}
  {\mathbb{E}} \; := \; \{ \,  | \vec{n} > := \exp( i \vec{n} \cdot \vec{\phi}) \; \; \vec{n}
  \in {\mathbb{Z}}^{2} \, \}
\end{equation}
given by:
\begin{equation}
  U_{T_{t}} | \vec{n} > \; = \; \exp( i \vec{n} \cdot \vec{\omega}
  t ) \, | \vec{n} >
\end{equation}
The geometrical phases in analysis emerges as the
\textbf{holonomy group $ HOL_{Stiefel} ( | n > ) $ in $ | \vec{n}
> $ }  of $ \omega_{Stiefel} $ where the
restriction to the vectors of the basis will be clarified in the
next section

\bigskip

\begin{example} \label{ex: Arnold cat}
\end{example}
ARNOLD CAT

Let us consider the automorphism of the torus  $ T^{2} \; := \;
{\mathbb{R}}^{2} / {\mathbb{Z}}^{2} $ identified by the linear
operator on $ {\mathbb{R}}^{2} $ whose matrix w.r.t. the
canonical basis is given by:
\begin{equation}
  C \; := \;  \begin{pmatrix}
    1 & 1 \\
    1 & 2 \
  \end{pmatrix}
\end{equation}
The passage to the Koopman representation involves  again the
unitary operator $ U_{C} $ over $ {\mathcal{H}} := L^{2} ( T^{2} ,
d \mu ( \phi_{1} , \phi_{2} ) ) $ identified by its action on the
basis:
\begin{equation}
  {\mathbb{E}} \; := \; \{ \,  | \vec{n} > := \exp( i \vec{n} \cdot \vec{\phi}) \; \; \vec{n}
  \in {\mathbb{Z}}^{2} \, \}
\end{equation}
given by:
\begin{equation}
  U_{C} | \vec{n} > \; = \; |  C \, \vec{n} >
\end{equation}
As in the example \ref{ex: Arnold cat}, the geometrical phases in
analysis emerges as the \textbf{holonomy group $ HOL_{Stiefel} (
| n > ) $ in $ | \vec{n}
> $ }  of $ \omega_{Stiefel} $ where the
restriction to the vectors of the basis will be clarified in the
next section.
\newpage
\section{The physical meaning of the classical geometric phase in term of moving frames}
In the previous section we have introduced, in the Koopman
representation of a classical dynamical system, the mathematical
setting underlying, in a quantum context, the Aharonov-Anandan
(non-adiabatic) quantum geometric phase.

In this paragraph we will explore the physical meaning of the
resulting (non adiabatic) classical geometric phase.

Let us start observing that, contrary to the quantum case, $
P{\mathcal{H}} $ doesn't represent the pure states of the system
but simply a particular class of physical observables observed up
to a phase \footnote{We are admitting here a class of physical
observables $ {\mathcal{H}  } $ broader than the usually accepted
one $ L^{\infty} (X, \mu) $} .

It is important, with this regard, to observe, that, contrary to
the quantum case, the neglected phase\emph{ do have here physical
meaning}, since $ | \psi > $ and $ \exp ( i \theta ) | \psi > $
are physically distinguishable.

Let us observe, furthermore, that the similarity between the
evolution equation of an \textbf{observable} $ | \psi > \in
{\mathcal{H}} $:
\begin{equation}
  | \psi ( t=1) > \; = \; U_{T}  | \psi ( t=0) >
\end{equation}
and the evolution-equation of Quantum Mechanics in the
\textbf{Schr\"{o}dinger picture} is deceptive since in our
classical Koopmanian situation we are in the \textbf{Heisenberg
picture} according to which the observables evolve with time
while the \textbf{state} $ \omega_{\mu} : L^{2} ( X \, , \, d \mu
) \mapsto {\mathbb{C}}$:
\begin{equation}
  \omega_{\mu} (f ) \; := \; \int_{X} f \, d \mu
\end{equation}
doesn't evolve with time.

Considered a basis of $ {\mathcal{H}} $:
\begin{equation}
   {\mathbb{E}} \; := \; \{ \,  | n > \, , \, n \in {\mathbb{N}}  \}
\end{equation}
let us consider the associated family of projectors

\begin{equation}
  P{\mathbb{E}} \; := \; \{ \,  | n > < n | \, , \, n \in {\mathbb{N}}  \}
\end{equation}
It determines the \textbf{frame} we use to  decompose and
consequentially analyze the particular class of  physical
observables constitued by the elements of $ {\mathcal{H}}$.

When such a \textbf{frame} is allowed \textbf{to move} for a
second following a loop $ \gamma \in C_{ |n > < n |}
(P{\mathcal{H}}) $:
\begin{equation}
  \gamma : [ 0 , 1 ] \mapsto P{\mathcal{H}}  \: : \: \gamma (t = 0)
   \, = \, \gamma (t = 1) \,  = \, | n > < n |
\end{equation}

the way the physical observables of $ {\mathcal{H}} $ observed,
through it, appear to us in the time interval  $( 0 , 1 ) $ is
the net result of two effects:
\begin{enumerate}
  \item  the real dynamical evolution of the frame
  \item  the motion of the frame we are performing making it to
  follow the prescribed loop
\end{enumerate}
At time t=1 the state of affairs of our observables aren't
identical to that we would have if we hadn't moved our frame: they
differ precisely by the geometric phase $
\tau_{\gamma}^{\omega_{Stiefel}} ( | n
> ) := \exp ( i \theta ) $ as it may be more clearly
visualized in the examples:

\begin{example} \label{ex:nonadiabatic geometric phase for the uncoupled oscillators}
\end{example}
NON ADIABATIC GEOMETRIC PHASE OF THE UNCOUPLED OSCILLATORS

At time $ t=1 $ we will have that:
\begin{equation}
  | \vec{n} (t =1 ) >  \; = \; e^{i \theta } U_{T_{t}} | \vec{n} > \; = \; \exp( i \vec{n} \cdot \vec{\omega}
  t \, + \, \theta ) \, | \vec{n} >
\end{equation}

\begin{example} \label{ex:nonadiabatic geometric phase for the Arnold cat}
\end{example}
ARNOLD CAT

At time $ t=1 $ we will have that:
\begin{equation}
  | \vec{n} (t=1)  >  \; = \; \exp ( i \theta ) U_{C} | \vec{n}  > \; = \; \exp ( i \theta ) |  C \, \vec{n}  >
\end{equation}
\newpage
\section{The recovering of Hannay angle in the adiabatic case}
The link between the fibre-bundle setting underlying the
Aharonov-Anandan (non-adiabatic) quantum geometric phase and the
fibre-bundle setting underlying Berry adiabatic quantum geometric
phase (discovered by Barry Simon) has been mathematically
characterized by A. Bohm amd Al\'{i} Mostafazadeh simply through a
suitable pull-back of the involved fibre bundle.

It is then reasonable that both the fibre-bundle setting
underlying Hannay adiabatic classical geometric phase and its
link with the non adiabatic classical geometric phase discusssed
in the previous sections may be simply obtained mimicking the
Bohm-Mostafazadeh-pullback.

Letc us consider an Hamiltonian integrable dynamical system
\cite{Arnold-Kozlov-Neishadt-93} with motion-equation, in the
angle-action variables $ ( \vec{I} \, , \, \vec{\phi} ) $, given
by:
\begin{equation}
  \frac{d \vec{I} }{dt} \; = \; 0
\end{equation}
\begin{equation}
  \frac{d \vec{\phi} }{dt} \; = \; \vec{\omega}(I)
\end{equation}
with:
\begin{equation}
 \vec{\omega} (I) \; = \; \frac{ \partial H }{ \partial \vec{I} }
\end{equation}

If the dynamical system in analysis has n degrees of freedoms,
and hence phase space 2-n dimensional, its Koopman description
involves the introduction of the Hilbert space:
\begin{equation}
  {\mathcal{H}} \; := \; L^{2} ( T^{n} , \frac{d \vec{\phi} }{(2
  \pi)^{n}})
\end{equation}
and of the Koopman unitary $ U_{t} $ identified by its action on
the basis:
\begin{equation}
  {\mathbb{E}} \; := \; \{ \,  | \vec{n} > := \exp( i \vec{n} \cdot \vec{\phi}) \; \; \vec{n}
  \in {\mathbb{Z}}^{n} \, \}
\end{equation}
given by:
\begin{equation}
  U_{t} | \vec{n} > \; = \; \exp( i \vec{n} \cdot \vec{\omega}
  t ) \, | \vec{n} >
\end{equation}

Let us suppose to alter the hamiltonian $ H \rightarrow H(
\vec{R} ) $ making the parameter $ \vec{R} $ to evolve
adiabatically  realizing a  loop in a suitable parameter space M.

In the adiabatic limit the basis:
\begin{equation}
  {\mathbb{E}}_{\vec{R}} \; := \; \{ \,  | \vec{n} , \vec{R} > \; \; \vec{n}
  \in {\mathbb{Z}}^{n} \, , \vec{R} \in M  \}
\end{equation}
continues to be formed by eigenvectors of $ U_{t} $.

Let us then introduce the following map $ f_{\vec{n}} :
{\mathbb{Z}}^{n} \times M \, \mapsto \, P  {\mathcal{H}} $:
\begin{equation}
  f_{\vec{n}} ( \vec{R} ) \; := \;  | \vec{n} , \vec{R} > < \vec{n} , \vec{R} |
\end{equation}
Obviously $  f_{\vec{n}} $ may also be seen as a  $
\vec{n}$-parametrized family of maps with domain M and codomain $
P{\mathcal{H}}$.

Let us then introduce the pullback-bundle $ f_{\vec{n}}^{\star}
{\mathcal{S}} ( {\mathcal{H}}  ) $ of the Stiefel bundle $
{\mathcal{S}} ( {\mathcal{H}}  ) ( P{\mathcal{H}} , U(1) ) $ by $
f_{\vec{n}} $ and let us denote by $ f ^{\star} \omega_{Stiefel}
$ the connection on the principal bundle $  f_{\vec{n}}^{\star}
{\mathcal{S}} ( {\mathcal{H}}  ) $ induced by the Stiefel
connection $ \omega_{Stiefel} $ through the pull-back operation.

Given a loop $ \gamma \in C_{\vec{R}=0} (M) $ the associated
Hannay phase is then $ \tau_{\gamma}^{f ^{\star}
\omega_{Stiefel}} $.

\end{document}